\documentclass[twocolumn]{autartcustom}    

\usepackage{graphicx}          
\usepackage{amsmath}
\usepackage{amsfonts,amssymb}
\usepackage[utf8]{inputenc}
\begin{document}

\begin{frontmatter}

\title{Output feedback stabilization of the linearized Korteweg-de Vries equation with right endpoint controllers\thanksref{footnoteinfo}} 

\thanks[footnoteinfo]{This research was funded by IZTECH BAP Grant 2017IYTE14.}

\author[Izm]{Ahmet Batal}\ead{ahmetbatal@iyte.edu.tr }\text{ and}    
\author[Izm]{Türker Özsarı}\ead{turkerozsari@iyte.edu.tr}              

\address[Izm]{Department of Mathematics, Izmir Institute of Technology, Urla, Izmir, 35430 Turkey }  


\begin{abstract}                          
In this paper, we prove the output feedback stabilization for the linearized Korteweg-de Vries (KdV) equation posed on a finite domain in the case the full state of the system cannot be measured.  We assume that there is a sensor at the left end point of the domain capable of measuring the first and second order boundary traces of the solution.  This allows us to design a suitable observer system whose states can be used for constructing boundary feedbacks acting at the right endpoint so that both the observer and the original plant become exponentially stable.  Stabilization of the original system is proved in the $L^2$-sense, while the convergence of the observer system to the original plant is also proved in higher order Sobolev norms.  The standard backstepping approach used to construct a left endpoint controller fails and presents mathematical challenges when building right endpoint controllers due to the overdetermined nature of the related kernel models.  In order to deal with this difficulty we use the method of \cite{BatalOzsari2018-1} which is based on using modified target systems involving extra trace terms.  In addition, we show that the number of controllers and boundary measurements can be reduced to one, with the cost of a slightly lower exponential rate of decay.  We provide numerical simulations illustrating the efficacy of our controllers.
\end{abstract}

\end{frontmatter}

\section{Introduction}
In this paper, we study the output feedback stabilization of the linearized Korteweg-de Vries (KdV) equation on a bounded domain $\Omega=(0,L)\subset \mathbb{R}$.  The linearized version of the model under consideration is given by
\begin{equation}\label{KdVBurgers}
 	\begin{cases}
 	\displaystyle u_{t} + u_{x} + u_{xxx} =0 \text { in } \Omega\times (0,T),\\
    u(0,t) = 0, u(L,t) = U(t), u_{x}(L,t)=V(t),\\
    u(x,0) = u_0(x) \text { in } \Omega,\\
 	\end{cases}
\end{equation} whereas the nonlinear version of this model is written with the main equation replaced by \begin{equation}\label{nonlinearKdV}u_t+u_x+u_{xxx}+uu_x=0.\end{equation}

In \eqref{nonlinearKdV}, $u=u(x,t)$ can for example model the evolution of the amplitude of a surface water wave in a finite length channel where energy to the system is put from the right end and left end of the system is free. The inputs $U(t)$ and $V(t)$ at the right end point of the boundary are feedback controllers to be constructed. The initial-boundary value problems \eqref{KdVBurgers} and \eqref{nonlinearKdV} with homogeneous boundary conditions ($U=V\equiv 0$) are both dissipative, since their solution satisfies $\frac{d}{dt}\|u(t)\|_{L^2(\Omega)}^2\le 0.$ However, this does not always guarantee exponential decay.  It is well-known that for some special domain lengths (so called \emph{critical lengths} for KdV) the solution does not need to decay to zero at all.  For example if $L=2\pi$, $u=1-\cos(x)$ is a (time independent) solution of \eqref{KdVBurgers} on $\Omega=(0,2\pi)$, but its $L^2-$norm is constant in $t$.  Therefore, introducing a stabilizing effect into the system is essential if one desires to steer the solution to zero.  See also \cite{Chu15} and \cite{Tang2016} for a detailed discussion of the relationship between stability and domain length.

If the state of the system can be measured at all times, one can attempt to construct exponentially stabilizing backstepping controllers for \eqref{KdVBurgers} and \eqref{nonlinearKdV}. A backstepping controller is generally constructed by using a transformation given by $w(x,t)=u(x,t)-\int_0^xk(x,y)u(y,t)dy,$ where $k$ is a kernel function which is chosen in such a way that the solution of \eqref{KdVBurgers} can be mapped to the solution of the following problem (so called ``\emph{target system}''):
\begin{equation*}
 	\begin{cases}
 	\displaystyle w_{t} + w_{x} + w_{xxx} + \lambda w = 0 & \text { in } \Omega\times \mathbb{R_+},\\
    w(0,t) = w(L,t) = w_{x}(L,t) = 0 & \text { in } \mathbb{R_+},\\
    w(x,0)=w_0(x)\equiv u_0-\int_0^xk(x,y)u_0(y)dy & \text { in } \Omega.
 	\end{cases}
\end{equation*}
The reason is that the solution of the above PDE model readily decays to zero, and if one can also show that the inverse of the backstepping transformation is bounded, then the decay of $w$ becomes equivalent to the decay of $u$.  Finding a suitable kernel which serves this purpose is the crucial step.  However, such an attempt to control from the right endpoint, when only one boundary condition is specified at the left, brings serious mathematical challenges since then the kernel is forced to satisfy the overdetermined PDE model given by
\begin{eqnarray} \label{kEqEEE}
  \nonumber k_{xxx} + k_{yyy} + k_{y} + k_{x}   = -\lambda k, y\in [0,x],x\in [0,L]\\
   k(x,x) = k(x,0)=k_y(x,0) = 0, \\
  \nonumber k_{x}(x,x) = \frac{\lambda}{3}x.
\end{eqnarray}
Unfortunately, the above PDE model does not have smooth solutions (see \cite{BatalOzsari2018-1} for a detailed discussion of this issue).  This problem is not present if one controls the system from the left endpoint \cite{Cerpa2013} or alternatively controls from the right with two boundary conditions specified at the left. The latter approach was used for instance in \cite{Tang2013} and \cite{Tang2015} where the controller acted from the right boundary condition while two (mixed type) boundary conditions were specified at the left.  However, usually the boundary conditions are determined by the intrinsic nature of the physical model, and one may not be able to choose the number of boundary conditions at a particular endpoint.  The novelty of the present article is that we are able to construct boundary feedback stabilizers acting from the opposite of the endpoint where only one boundary condition is specified.

Two approaches were proposed in order to overcome the difficulty associated with the overdetermined kernel model. Coron \& Lü \cite{cor14} replaced \eqref{kEqEEE} posed on a triangle with an equivalent PDE model posed on the rectangle $[0,L]\times [0,L]$ and showed that this kernel PDE model has a rough ($H^1$) solution.  However, their result relies on the exact controllability of the linear KdV equation, which does not hold on domains of\emph{ critical }lengths.  They managed to get high decay rates for domains of \emph{uncritical} lengths.  The second approach due to \cite{BatalOzsari2018-1} is a direct method which does not rely on any controllability result. It is based on constructing a backstepping controller which uses a modified kernel model disregarding one of the boundary conditions in \eqref{kEqEEE}:
\begin{eqnarray} \label{kEq}
\nonumber k_{xxx} + k_{yyy} + k_{y} + k_{x}   = -\lambda k, y\in [0,x],x\in [0,L]\\
   k(x,x) = k(x,0) = 0, \\
 \nonumber k_{x}(x,x) = \frac{\lambda}{3}x.
\end{eqnarray} In \cite{BatalOzsari2018-1} it is proven that the exponential stability can still be achieved by using such a kernel with the cost of a low exponential rate of decay.  The slower decay is due to the fact that disregarding a boundary condition from \eqref{kEqEEE} changes the target system in such a way that its main equation involves a trace term which depends on the kernel. Although this trace term badly affects the decay, its effect can be eliminated by choosing $\lambda$ sufficiently small in which case one can still obtain an exponential decay but not with an arbitrarily large rate. For more details see \cite[Section 2.1]{BatalOzsari2018-1}. This approach has the advantage that it is independent of whether the domain length is critical or not.  The existence as well as the smoothness of the kernel $k$ satisfying \eqref{kEq} was previously proved in \cite[Lemma 2.1]{BatalOzsari2018-1}:
\begin{lem}[\cite{BatalOzsari2018-1}]
There exists a $C^\infty$-function $k$ that solves the boundary value problem \eqref{kEq}.
\end{lem}
The proof of the above lemma was done in \cite{BatalOzsari2018-1} in two steps. The first step was to show that $k$ solves \eqref{kEq} if and only if $G=G(s,t)$ solves the integral equation
\begin{multline}\label{GepsInt1}
  {G}(s,t) = \frac{{\lambda}}{3}st\\
  +\frac{1}{3}\int_0^t\int_0^s\int_0^\omega (-{G}_{ttt} + 3{G}_{stt}  - {G}_{t} -\lambda {G})(\xi,\eta)d\xi d\omega d\eta,
\end{multline} where  $t\equiv y$, $s\equiv x-y$, and $G(s,t)\equiv k(x,y)$. The second step was to obtain the smooth solution of \eqref{GepsInt1}  via a successive approximation technique and uniform boundedness analysis of the subsequent series.

In both \cite{cor14} and \cite{BatalOzsari2018-1}, it was assumed that the state of the system could be measured at all times.  Unfortunately, this is not always the case.  For instance, if one has no access to the medium, a controller that requires measurement of the full state of the original system may not be constructed. In such a case, one generally first constructs an observer system that estimates the plant if some partial information such as a boundary measurement is available.  The advantage is that the observer can be controlled since its full state is available unlike the original plant.  This implies that the original plant can be stabilized by the same controller applied to the boundary of the observer.  From the mathematical point of view, the question is the following.
\begin{prob}\label{mainprob}
Can you write a boundary feedback system with exponential stability, say with the unknown $\hat{u}$, such that this system (observer) estimates the solution of the original plant with the same controller which uses the states of the observer?
\end{prob}
In this paper, we will assume that there are sensors at the left end point of the channel capable of measuring the boundary traces $u_x(0,t)$ and $u_{xx}(0,t)$.
In order to answer Problem \ref{mainprob}, we introduce and stabilize the following \emph{observer} system whose boundary feedback will also be applied to the original plant:
\begin{equation}\label{observer} \left\{ \begin{array}{ll}
        \hat{u}_t+\hat{u}_x+\hat{u}_{xxx}+P_1(x)\left(u_x(0,t)-\hat{u}_x(0,t)\right)\\
        +P_2(x)\left( u_{xx}(0,t)-\hat{u}_{xx}(0,t)\right)=0,\text { in } \Omega\times (0,T),\\
        \hat{u}(0,t)=0,\,\hat{u}(L,t)={U}(t),\,\hat{u}_x(L,t)={V}(t),\text { in } (0,T),\\
        \hat{u}(x,0)=\hat{u}_0(x),\text { in } \Omega.\end{array} \right.
\end{equation}
Note that the error $\tilde{u}=\hat{u}-u$ satisfies the PDE model given by
\begin{equation}\label{error} \left\{ \begin{array}{ll}
        \tilde{u}_t+\tilde{u}_x+\tilde{u}_{xxx}\\
        =P_1(x)\tilde{u}_x(0,t)+P_2(x)\tilde{u}_{xx}(0,t)\text { in } \Omega\times (0,T);\\
        \tilde{u}(0,t)=0, \tilde{u}(L,t)=0,\tilde{u}_x(L,t)=0 \text { in } (0,T);\\
        \tilde{u}(x,0)=u_0(x)-\hat{u}_0(x) \text { in } \Omega. \end{array} \right.
\end{equation}
$P_1(x)$ and $P_2(x)$  are observer gains in \eqref{observer} and \eqref{error}, which are chosen in such a way that the solution $\hat{u}$ of the estimator can be later controlled and moreover the error is enforced to go to zero as $t$ gets larger (see Section \ref{ObserverDes}).  In some sense, we want to control the error, too.  This is achieved by using a bounded invertible (backstepping) transformation  in the form
\begin{equation}\label{transtildew}
\tilde{u}(x,t)=\tilde{w}(x,t)-\int_0^x p(x,y)\tilde{w}(y,t)dy
\end{equation}
by mapping the error system to the (exponentially stable) target system given by
\begin{equation}\label{tildew} \left\{ \begin{array}{ll}
        \tilde{w}_t+\tilde{w}_x+\tilde{w}_{xxx}+\tilde{\lambda} \tilde{w}=0, \text { in } \Omega\times (0,T),\\
        \tilde{w}(0,t)=0,\,\tilde{w}(L,t)=0,\,\\\tilde{w}_x(L,t)=\int_0^Lp_x(L,y)\tilde{w}(y,t)dy, \text { in } (0,T),\\
        \tilde{w}(x,0)=\tilde{w}_0(x), \text { in } \Omega \end{array} \right.
\end{equation}
where $\tilde{\lambda}>0$. Computing the relevant partial derivatives of both sides of \eqref{transtildew}, applying integration by parts and using the given boundary conditions, it can be shown that the desired target system \eqref{tildew} is obtained if $P_1(x):=p_y(x,0),\; P_2(x):=-p(x,0)$ and $p(x,y)$ satisfies the following PDE model on $\Delta$:
\begin{equation}\label{p} \left\{ \begin{array}{ll}
        p_{xxx}+p_{yyy}+p_x+p_y=\tilde{\lambda} p,\\
        p(L,y)=0, \; p(x,x)=0, \\
        p_x(x,x)=-\frac{\tilde{\lambda}}{3}(x-L).\end{array} \right.
\end{equation}
Existence of a solution to \eqref{p} as well as the exponential decay of \eqref{tildew} are shown in Section \ref{ObserverDes} below.
\subsection{A few more words on the literature}
Recently, \cite{Marx18} proved the output feedback stabilization of the Korteweg-de Vries equation subject to the boundary conditions $u(0,t)=U(t),\,u_x(L,t)=u_{xx}(L,t)=0$ by using the partial measurement $y(t)=u(L,t)$ .  Here the left end boundary input $U(t)$ is a controller (stabilizer) obtained by using the backstepping method.  This controller uses only the state values of the observer. Prior to this work, the same authors \cite{Marx14} proved the output feedback stabilization of the Korteweg-de Vries equation subject to the boundary conditions $u(0,t)=U(t),\,u(L,t)=u_{x}(L,t)=0$ by using the partial measurement $y(t)=u_{xx}(L,t)$. The same problem in the nonlinear case was studied by \cite{Hasan2016}.  We should also mention some important work related to the control and stabilization of the KdV equation.  Exact boundary controllability of the linear and nonlinear KdV equations with the same type of boundary conditions as in \eqref{KdVBurgers} was studied by \cite{Rosier1997}, \cite{Cor2004}, \cite{Zhang99}, \cite{Glass08}, \cite{Cerpa07}, \cite{Cerpa09}, \cite{RosZha09}, and \cite{Glass10}.  Stabilization of solutions of the KdV equation with a localised interior damping was achieved by \cite{Perla2002}, \cite{Pazo05}, \cite{Mass07}, and \cite{Balogh2000}. There are also some results achieving stabilization of the KdV equation by using predetermined local boundary feedbacks, see for instance \cite{Liu2002} and \cite{Jia2016}.
\subsection{Preliminaries, notation, and main result}
Before we state our main results, let us give some important facts and notations that will be needed later. To this end, let $\eta$ be a $C^\infty$-function and $\Upsilon_{\eta}:H^l(\Omega)\rightarrow H^l(\Omega)$ ($l\ge 0$) be the integral operator defined by $(\Upsilon_{\eta}\varphi)(x):=\int_0^x{\eta}(x,y)\varphi(y)dy,$ where $H^l(\Omega)$ denotes the $L^2-$based Sobolev spaces with $H^0(\Omega)=L^2(\Omega)$. Then the following result holds true \cite{Liu03}, \cite{BatalOzsari2018-1}:
\begin{lem}\label{inverselem}
 $I-\Upsilon_{\eta}$ is invertible with a bounded inverse from $H^l(\Omega)\rightarrow H^l(\Omega)$ ($l\ge 0$). Moreover,  $(I-\Upsilon_{\eta})^{-1}$ can be written as $I+\Phi$, where $\Phi$ is a bounded operator from $L^2(\Omega)$ into $H^l(\Omega)$ for $l=0,1,2$ and from $H^{l-2}(\Omega)$ into $H^{l}(\Omega)$ for $l> 2$.
\end{lem}
For a given function $\varphi$, we say it satisfies the (higher order) compatibility conditions (see e.g., \cite[Definition 1.1]{BSZ03}) if\begin{equation}\label{compa}
                                  \varphi(\bar{x})=\varphi'''(\bar{x})+\varphi'(\bar{x})=0, \bar{x}=0,L.
                                \end{equation}
We also set $X_T^s=C([0,T];H^s(\Omega))\cap L^2(0,T;H^{s+1}(\Omega))$ for representing solution spaces for $s\ge 0$.
In what follows, we will write $ A\lesssim B$ to denote an inequality $A\leq cB$ where $c>0$ may only depend on the fixed parameters of the problem under consideration which are not of interest. The main result of the paper is stated in the following theorem:
\begin{thm}\label{mainthm}
  Let $T>0$, $u_0, \hat{u}_0\in H^{6}(\Omega)$ with $u_0(0)=u_0(L)=0$, $p$ and $k$ be the smooth kernels solving \eqref{p} and \eqref{kEq}, respectively.  Let also $(I-\Upsilon_p)^{-1}\tilde{u}_0=\tilde{w}_0$ satisfy the compatibility conditions \eqref{compa}. Then, the plant-observer-error (POE) system given in \eqref{KdVBurgers}, \eqref{observer}, \eqref{error} has a solution $(u,\hat{u},\tilde{u})\in X_T^3\times X_T^3\times X_T^6$ with right endpoint boundary controllers
  $$U(t) := [\Upsilon_k\hat{u}](L,t)\text{ and }V(t): =[\Upsilon_{k_x}\hat{u}](L,t).$$
 Moreover, there exist $\alpha>\kappa>0$ such that the decay rate estimates
\begin{eqnarray}
\nonumber\|u(t)\|_{L^2(\Omega)}\lesssim \left(\|\hat{u}_0\|_{L^2(\Omega)}+\|u_0-\hat{u}_0\|_{H^3(\Omega)}\right)e^{-\kappa t}\\
+\|u_0-\hat{u}_0\|_{L^2(\Omega)}e^{-\alpha t},\\
\|\hat{u}(t)\|_{L^2(\Omega)} \lesssim \left(\|\hat{u}_0\|_{L^2(\Omega)}+\|u_0-\hat{u}_0\|_{H^3(\Omega)}\right)e^{-{\kappa} t},\\
\|u(t)-\hat{u}(t)\|_{L^2(\Omega)} \lesssim \|u_0-\hat{u}_0\|_{L^2(\Omega)}e^{-{\alpha}t},\\
\|u(t)-\hat{u}(t)\|_{H^3(\Omega)} \lesssim \|u_0-\hat{u}_0\|_{H^3(\Omega)}e^{-{\alpha}t}
\end{eqnarray} hold true for $t\in [0,T].$
\end{thm}

\section{Linearized model}
\subsection{Wellposedness} The initial step is to prove the wellposedness of the target error system \eqref{tildew}.  To this end, we first consider the following open loop system instead of \eqref{tildew} for a moment:
\begin{equation}\label{q-prob} \left\{ \begin{array}{ll}
        {\tilde w}_t+{\tilde w}_x+{\tilde w}_{xxx}+\tilde{\lambda} {\tilde w}=0, \text { in } \Omega\times (0,T),\\
        {\tilde w}(0,t)=0,\,{\tilde w}(L,t)=0,\,{\tilde w}_x(L,t)=h(t),\\
        {\tilde w}(x,0)={\tilde w}_0(x), \text { in } \Omega, \end{array} \right.
\end{equation} where $h\in H^1(0,T)$, ${\tilde w}_0\in H^3(\Omega)$ satisfy the compatibility conditions ${\tilde w}_0(0)=0, {\tilde w}_0(L)=0$.  The well-posedness of \eqref{q-prob} was obtained in \cite[Lemma 3.3]{BSZ03}, and one has ${\tilde w}\in X_T^3$ together with ${\tilde w}_t\in X_T^0$.
\begin{lem}[\cite{BSZ03}]\label{BSZ03lem} For given $T>0$, let $h\in H^1(0,T)$, ${\tilde w}_0\in H^3(\Omega)$ satisfy the compatibility conditions ${\tilde w}_0(0)=0, {\tilde w}_0(L)=0$.  Then equation \eqref{q-prob} has a unique solution ${\tilde w}$ in $X_T^3$ with ${\tilde w}_t\in X_T^0$ such that the following estimates hold true:
 \begin{equation*}
\left(\|{\tilde w}\|_{X_T^3}+\|{\tilde w}_t\|_{{X_T^0}}\right) \le C \left(\|{\tilde w}_0\|_{H^3(\Omega)}+\|h\|_{H^1(0,T)}\right).
\end{equation*}
\end{lem}
Note that in \eqref{tildew}, the boundary condition $\tilde{w}_x(L,t)=\int_0^Lp_x(L,y){\tilde w}(y,t)dy$ is of feedback type.  This corresponds to a closed loop version of \eqref{q-prob} where $h(t)=h({\tilde w})(t)=\int_0^Lp_x(L,y){\tilde w}(y,t)dy$. The wellposedness of the closed loop problem will be treated by using a fixed point argument.  To achieve this, we define the Banach space
$Q_T\equiv \{{\tilde w}\in X_T^3\,|\,{\tilde w}_t\in X_T^0\}$ and its complete metric subspace $\tilde{Q}_T=\{{\tilde w}\in Q_T\,|\,{\tilde w}(\cdot,0)={\tilde w}_0(\cdot)\}$ with the metric induced from the norm of $Q_T.$ Observe that given ${\tilde w}^*\in \tilde{Q}_T$, since $p$ is a smooth solution of \eqref{p}, one has $h({\tilde w}^*)(\cdot)=\int_0^Lp_x(L,y){\tilde w}^*(y,\cdot)dy\in H^1(0,T)$. Indeed,
\begin{multline}\|h({\tilde w}^*)\|_{H^1(0,T)}=\left\|\int_0^Lp_x(L,y){\tilde w}^*(y,\cdot)dy\right\|_{H^1(0,T)}\\
\le \sqrt{T}\|p_x(L,\cdot)\|_{L^2(\Omega)}\left(\|{\tilde w}^*\|_{X_T^0}+\|{\tilde w}_t^*\|_{X_T^0}\right)<\infty. \end{multline}
Now, we  replace the boundary condition ${\tilde w}_x(L,t)=h({\tilde w})(t)$ with ${\tilde w}_x(L,t)=h({\tilde w}^*)(t)$ for fixed ${\tilde w}^*\in  X_T^3$. This is nothing but the problem given in \eqref{q-prob} which has a unique solution by Lemma \ref{BSZ03lem}.
This defines an operator $\Gamma: \tilde{Q}_T\rightarrow \tilde{Q}_T$ given by $\Gamma({\tilde w}^*)={\tilde w}$. Regarding the closed loop problem \eqref{tildew}, it is now enough to show that $\Gamma$ has a fixed point. Let ${\tilde w}_1,{\tilde w}_2\in \tilde{Q}_T$. Using the estimate in Lemma \ref{BSZ03lem}, we have
\begin{multline}d(\Gamma({\tilde w}_1),\Gamma({\tilde w}_2))_{\tilde{Q}_T}=\|\Gamma({\tilde w}_1)-\Gamma({\tilde w}_2)\|_{Q_T}\\
\le C\|h({\tilde w}_1)(\cdot)-h({\tilde w}_2)(\cdot)\|_{H^1(0,T)}\\
\le CT\|{\tilde w}_1-{\tilde w}_2\|_{Q_T}=CTd({\tilde w}_1,{\tilde w}_2)_{\tilde{Q}_T}.
\end{multline}  Note that by choosing $T$ sufficiently small we can make the constant at the right hand side of the above inequality less than 1.  Now, unleashing the Banach fixed point theorem, we obtain the existence of a unique local solution ${\tilde w}\in Q_T$. This implies the local well-posedness for the target error system \eqref{tildew}.
In order to show that the local solution is indeed global, it is enough to prove that the local solution stays uniformly bounded in time. But this readily follows from the stabilization estimates given in Section \ref{ObserverDes} below.  Now, by using the transformation in \eqref{transtildew}, we obtain the wellposedness of the error system \eqref{error}.

We prove in Lemma \ref{tracelem} below that $\tilde{w}_{x}(0,\cdot), \tilde{w}_{xx}(0,\cdot)\in L^2(0,T).$ Moreover,
for $\tilde{w}_0\in H^6(\Omega)$, we have $z_0:=-\tilde{w}_0'-\tilde{w}_0'''-\tilde{\lambda}\tilde{w}_0\in H^3(\Omega)$ satisfying the compatibility conditions. Introducing $z=\tilde{w}_t$, we observe that $z$ satisfies the main equation as well as the boundary conditions of \eqref{tildew} but with initial condition $z(x,0)=z_0$.  Applying the above arguments to $z$, we deduce that $\tilde{w}_t=z\in X_T^3$.  Moreover, we have $\tilde{w}_{xt}(0,t)=z_x(0,t), \tilde{w}_{xxt}(0,\cdot)=z_{xx}(0,\cdot)\in L^2(0,T)$.  Therefore, the right hand side of \eqref{hatw} can be written as $a(x)\hat{w}_x(0,t)+f(x,t)$ with $a(x)=k_y(x,0)$ and $f(x,t)=-\Psi_1(x)\tilde{w}_x(0,t)-\Psi_2(x)\tilde{w}_{xx}(0,t)$ such that $f\in W^{1,2}(0,T;H^\infty(\Omega))$.  Well-posedness of this problem was studied in \cite[Lemma 3.3]{BSZ03}, and for given $\hat{w}_0\in H^3(\Omega)$ satisfying the compatibility, one has  $\hat{w}\in X_T^3.$ Now, by the invertibility of \eqref{transhatw} due to Lemma \ref{inverselem}, we obtain the wellposedness of the observer system \eqref{observer} so that $\hat{u}\in X_T^3.$  Combining the wellposedness of \eqref{error} and \eqref{observer}, we obtain the wellposedness of the original system and conclude that ${u}=X_T^3.$

\subsection{Stabilization}\label{ObserverDes}

Note that with the change of variables $\tilde{x}\equiv L-y$ and $\tilde{y}\equiv L-x$  and $k(\tilde{x},\tilde{y})=p(x,y)$, it is easy to see that $k$ is the $C^\infty$ kernel which solves \eqref{kEq}, where $x$, $y$, and $\lambda$, replaced by $\tilde{x}$, $\tilde{y},$ and $\tilde{\lambda}.$
Note also that by this transformation we see that $p_x(L,y)=-k_{\tilde{y}}(\tilde{x},0)$, and in \cite[Lemma 2.5]{BatalOzsari2018-1} it is shown that for suitably small, $\tilde{\lambda}>0$, the quantity $\tilde{\lambda}-\frac{1}{2}\|k_{\tilde{y}}(\cdot,0)\|_{L^2(\Omega)}^2$ is strictly greater than zero. Therefore choosing $\tilde{\lambda}$ sufficiently small, we can guarantee that $\alpha \equiv \tilde{\lambda}-\frac{1}{2}\|p_x(L,\cdot )\|^2_{L^2(\Omega)}>0$. We need the following lemma:
\begin{lem}\label{tracelem}
Let $\tilde{w}$ be the solution of \eqref{tildew}. Then the following inequalities hold:
\begin{eqnarray}
\label{ineq1}\|\tilde{w}\|_{L^2(\Omega)} \leq \|\tilde{w}_0\|_{L^2(\Omega)}e^{-\alpha t}, \\
\label{ineq2}|\tilde{w}_x(0,t)|+|\tilde{w}_{xx}(0,t)|+\|\tilde{w}\|_{H^3(\Omega)}\lesssim \|\tilde{w}_0\|_{H^3(\Omega)}e^{-\alpha t}.
\end{eqnarray}
\end{lem}
We multiply \eqref{tildew} by $\tilde{w}$ and integrate over $\Omega$. Applying integration by parts and boundary conditions we obtain
\begin{multline*}
\frac{1}{2}\frac{d}{dt}\|\tilde{w}(t)\|_{L^2(\Omega)}^2+\tilde{\lambda}\|\tilde{w}(t)\|_{L^2(\Omega)}^2+\frac{1}{2}
|\tilde{w}_x(0,t)|^2\\
=\frac{1}{2}|\tilde{w}_x(L,t)|^2,
\end{multline*}
which, together with \eqref{tildew}, implies
\begin{multline*}
\frac{1}{2}\frac{d}{dt}\|\tilde{w}(t)\|_{L^2(\Omega)}^2+\tilde{\lambda}\|\tilde{w}(t)\|_{L^2(\Omega)}^2\\
\leq \frac{1}{2}\bigg(\int_0^Lp_x(L,y)\tilde{w}(y,t)dy\bigg)^2.
\end{multline*}
Applying the Cauchy-Schwarz inequality to the right hand side we see that
\begin{equation*}
\frac{1}{2}\frac{d}{dt}\|\tilde{w}(t)\|_{L^2(\Omega)}^2+\big(\tilde{\lambda}-\frac{1}{2}\|p_x(L,\cdot )\|^2_{L^2(\Omega)}\big)\|\tilde{w}(t)\|_{L^2(\Omega)}^2\leq 0,
\end{equation*}
which gives \eqref{ineq1}.

In order to prove \eqref{ineq2}, we first differentiate \eqref{tildew} with respect to $t$, then multiply by $\tilde{w}_{t}$ and integrate over $\Omega$. Using integration by parts and boundary conditions as well, we see that
\begin{multline}
\frac{1}{2}\frac{d}{dt}\|\tilde{w}_t(t)\|_{L^2(\Omega)}^2+\tilde{\lambda}\|\tilde{w}_t(t)\|_{L^2(\Omega)}^2+\frac{1}{2}
|\tilde{w}_{tx}(0,t)|^2\\
=\frac{1}{2}|\tilde{w}_{tx}(L,t)|^2.
\end{multline}
Moreover by \eqref{tildew} we have $\tilde{w}_{tx}(L,t)=\int_0^Lp_x(L,y)\tilde{w}_t(y,t)dy$. Hence we obtain
\begin{multline*}
\frac{1}{2}\frac{d}{dt}\|\tilde{w}_t(t)\|_{L^2(\Omega)}^2+\tilde{\lambda}\|\tilde{w}_t(t)\|_{L^2(\Omega)}^2\\
\leq \frac{1}{2}\bigg(\int_0^Lp_x(L,y)\tilde{w}_t(y,t)dy\bigg)^2.
\end{multline*}
Applying the Cauchy-Schwarz inequality to the right hand side we get
\begin{equation}
\frac{1}{2}\frac{d}{dt}\|\tilde{w}_t(t)\|_{L^2(\Omega)}^2+\big(\tilde{\lambda}-\frac{1}{2}\|p_x(L,\cdot )\|^2_{L^2(\Omega)}\big)\|\tilde{w}_t(t)\|_{L^2(\Omega)}^2\leq 0,
\end{equation}
which implies
\begin{equation}\label{wtdecay}
\|\tilde{w}_t(t)\|_{L^2(\Omega)}\leq \|\tilde{w}_t(0)\|_{L^2(\Omega)}e^{-\alpha t}\leq \|\tilde{w}_0\|_{H^3(\Omega)}e^{-\alpha t}
\end{equation}
 since $ \|\tilde{w}_t(0)\|_{L^2(\Omega)}=\|\tilde{w}_0'+\tilde{w}_0'''+\tilde{\lambda}\tilde{w}_0\|\leq \|\tilde{w}_0\|_{H^3(\Omega)}.$
 On the other hand, by \eqref{tildew} we also have
 \begin{multline}\label{wxxxest}
 \|\tilde{w}_{xxx}(t)\|_{L^2(\Omega)}^2\\
 \leq 3\big(\|\tilde{w}_{x}(t)\|_{L^2(\Omega)}^2+\tilde{\lambda}\|\tilde{w}(t)\|_{L^2(\Omega)}^2+\|\tilde{w}_{t}(t)\|_{L^2(\Omega)}^2\big).
 \end{multline}
 Applying $\epsilon$-Young's inequality to the square of the right hand side of the Gagliardo-Nirenberg inequality
 \begin{equation*}
 \|\tilde{w}_{x}(t)\|_{L^2(\Omega)}\leq \|\tilde{w}_{xxx}(t)\|^{\frac{1}{3}}_{L^2(\Omega)}\|\tilde{w}(t)\|^{\frac{2}{3}}_{L^2(\Omega)},
 \end{equation*}
 we also obtain
 \begin{equation}\label{wxest}
 \|\tilde{w}_{x}(t)\|_{L^2(\Omega)}^2\leq \epsilon \|\tilde{w}_{xxx}(t)\|_{L^2(\Omega)}^2+c_\epsilon\|\tilde{w}(t)\|_{L^2(\Omega)}^2
 \end{equation}
 for $\epsilon>0$.
Combining \eqref{wxxxest} and \eqref{wxest}, and choosing $\epsilon$ small enough, we see that
\begin{equation}\label{anewold1}
\|\tilde{w}_{xxx}(t)\|_{L^2(\Omega)} \lesssim \|\tilde{w}(t)\|_{L^2(\Omega)}+\|\tilde{w}_{t}(t)\|_{L^2(\Omega)}.
\end{equation}
Hence
\begin{equation}\label{anewold2}\|\tilde{w}(t)\|_{H^3(\Omega)}\lesssim \|\tilde{w}(t)\|_{L^2(\Omega)}+\|\tilde{w}_{t}(t)\|_{L^2(\Omega)},\end{equation}
which, together with \eqref{ineq1} and \eqref{wtdecay}, implies
\begin{equation}\label{part2}
\|\tilde{w}\|_{H^3(\Omega)} \lesssim \|\tilde{w}_0\|_{H^3(\Omega)}e^{-\alpha t}.
\end{equation}
To obtain the second part of inequality \eqref{ineq2}, we multiply \eqref{tildew} by $(L-x)\tilde{w}_{xx}$ and integrate over $\Omega$. Applying integration by parts and boundary conditions, we obtain
\begin{multline*}
\tilde{w}_x^2(0,t)+\tilde{w}_{xx}^2(0,t)\\
=\frac{2}{L}\int_0^L \big( (L-x)\tilde{w}_t\tilde{w}_{xx}+\frac{1}{2}\tilde{w}_{x}^2+\frac{1}{2}\tilde{w}_{xx}^2+\tilde{\lambda}(L-x)\tilde{w}\tilde{w}_{xx}\big) dx.
\end{multline*}
Using Cauchy-Schwarz and Young's inequalities on the first and last term of the right hand side, we see that
\begin{equation}
|\tilde{w}_x(0,t)|^2+|\tilde{w}_{xx}(0,t)|^2 \lesssim \|\tilde{w}_t\|_{L^2(\Omega)}^2+\|\tilde{w}\|_{H^3(\Omega)}^2,
\end{equation}
which, together with \eqref{wtdecay} and \eqref{part2}, implies \eqref{ineq2}.

Now for $\hat{u}$, we apply the backstepping transformation
\begin{equation}\label{transhatw}
\hat{w}=\hat{u}-\int_0^x k(x,y)\hat{u}(y,t)dy
\end{equation}
where $k$ is the kernel \cite[Lemma 2.1]{BatalOzsari2018-1} which solves \eqref{kEq}.
Choosing \begin{equation}\label{controllers}\begin{array}{ll}
{U}(t)=\int_0^L k(L,y)\hat{u}(y,t)dy,\\
{V}(t)=\int_0^L k_x(L,y)\hat{u}(y,t)dy. \end{array}
\end{equation} in \eqref{observer}, we obtain the following equation for $\hat{w}$:
\begin{equation}\label{hatw} \left\{ \begin{array}{ll}
        \hat{w}_t+\hat{w}_x+\hat{w}_{xxx}+\lambda \hat{w}\\
        =k_y(x,0)\hat{w}_x(0,t)-\Psi_1(x)\tilde{w}_x(0,t)-\Psi_2(x)\tilde{w}_{xx}(0,t),\\
        \hat{w}(0,t)=0, \; \hat{w}(L,t)=0, \; \hat{w}_x(L,t)=0, \end{array} \right.
\end{equation}
where $\Psi_i(x)\equiv P_i(x)-\int_0^x P_i(y)k(x,y)dy$ for $i\in \{1,2\}$. Multiplying \eqref{hatw} by $\hat{w}$ and integrating over $\Omega$, we obtain
\begin{multline*}
\frac{1}{2}\frac{d}{dt}\|\hat{w}(t)\|_{L^2(\Omega)}^2+\lambda\|\hat{w}(t)\|_{L^2(\Omega)}^2+\frac{1}{2}
|\hat{w}_x(0,t)|^2\\
=\hat{w}_x(0,t)\int_0^L k_y(x,0)\hat{w}(x,t)dx\\
-\tilde{w}_x(0,t)\int_0^L\Psi_1(x) \hat{w}(x,t)dx-\tilde{w}_{xx}(0,t)\int_0^L\Psi_2(x) \hat{w}(x,t)dx.
\end{multline*}
 Applying $\epsilon$-Young's and Cauchy-Schwarz inequalities to the right hand side, for any $\epsilon$ we get
\begin{equation}\label{estenergy}
\frac{1}{2}\frac{d}{dt}\|\hat{w}(t)\|_{L^2(\Omega)}^2+\kappa \|\hat{w}(t)\|_{L^2(\Omega)}^2
\leq \frac{1}{2\epsilon}\left[\tilde{w}^2_x(0,t)+\tilde{w}^2_{xx}(0,t)\right],
\end{equation}
where $$\kappa\equiv \lambda-\frac{1}{2}\|k_y(\cdot,0)\|_{L^2(\Omega)}^2-
\frac{1}{2}\epsilon\big(\|\Psi_1\|_{L^2(\Omega)}^2 + \|\Psi_2\|_{L^2(\Omega)}^2 \big).$$ By \cite[Lemma 2.5]{BatalOzsari2018-1} we know that for sufficiently small $\lambda$, the quantity $\lambda-\frac{1}{2}\|k_y(\cdot,0)\|_{L^2(\Omega)}^2>0$. Therefore choosing $\epsilon$ sufficiently small we can make the coefficient $\kappa> 0$. Moreover, since $\alpha \equiv \tilde{\lambda}-\frac{1}{2}\|p_x(L,\cdot )\|^2_{L^2(\Omega)}$ and $p_x(L,y)=k_y(x,0)$ , choosing $\lambda=\tilde{\lambda}$ if necessary, we can assume $\alpha>\kappa$. Inequality \eqref{estenergy} and \eqref{ineq2} imply
\begin{equation*}
\frac{1}{2}\frac{d}{dt}\|\hat{w}(t)\|_{L^2(\Omega)}^2+\kappa \|\hat{w}(t)\|_{L^2(\Omega)}^2\lesssim \|\tilde{w}_0\|^2_{H^3(\Omega)}e^{-2\alpha t}.
\end{equation*}
Using the assumption $\alpha>\kappa$, multiplying both sides of the above inequality by $e^{2\kappa t}$ and taking the integral of both sides from $0$ to $t$ we can easily see that $$\|\hat{w}(t)\|_{L^2(\Omega)}^2\lesssim \left(\|\hat{w}_0\|_{L^2(\Omega)}^2 +\|\tilde{w}_0\|^2_{H^3(\Omega)}\right)e^{-2\kappa t},$$ which is equivalent to saying
\begin{equation}\label{expdecay1}
\|\hat{w}(t)\|_{L^2(\Omega)}\lesssim \left(\|\hat{w}_0\|_{L^2(\Omega)}+\|\tilde{w}_0\|_{H^3(\Omega)}\right) e^{-\kappa t}.
\end{equation}

On the other hand, both of the transformations given in \eqref{transtildew} and \eqref{transhatw} are bounded with bounded inverses  by Lemma \ref{inverselem}. Therefore, we have
\begin{equation}\label{normeq}\begin{array}{ll}
\|\tilde{u}(t)\|_{H^3(\Omega)} \lesssim \|\tilde{w}(t)\|_{H^3(\Omega)}, & \|\tilde{w}_0\|_{H^3(\Omega)} \lesssim \|\tilde{u}_0\|_{H^3(\Omega)},\\
\|\hat{u}(t)\|_{L^2(\Omega)} \lesssim \|\hat{w}(t)\|_{L^2(\Omega)}, & \|\hat{w}_0\|_{L^2(\Omega)} \lesssim \|\hat{u}_0\|_{L^2(\Omega)}.
\end{array}
\end{equation}
Combining \eqref{expdecay1} and \eqref{normeq}, we achieve
\begin{equation}\label{uhat1}
\|\hat{u}\|_{L^2(\Omega)}\lesssim \left(\|\hat{u}_0\|_{L^2(\Omega)}+\|u_0-\hat{u}_0\|_{H^3(\Omega)}\right)e^{-\kappa t}.
\end{equation}
Moreover,  \eqref{ineq1}, \eqref{ineq2} and \eqref{normeq} also imply
\begin{eqnarray}\label{utilde1}
\|u-\hat{u}\|_{L^2(\Omega)} &\lesssim& \|u_0-\hat{u}_0\|_{L^2(\Omega)}e^{-\alpha t},\\
\|u-\hat{u}\|_{H^3(\Omega)} &\lesssim& \|u_0-\hat{u}_0\|_{H^3(\Omega)}e^{-\alpha t}.
\end{eqnarray}
Using \eqref{uhat1}-\eqref{utilde1} together with the triangle inequality, we obtain
\begin{multline}\|{u}\|_{L^2(\Omega)}=\|\hat{u}+\tilde{u}\|_{L^2(\Omega)}\le \|\hat{u}\|_{L^2(\Omega)}+\|u-\hat{u}\|_{L^2(\Omega)}\\
\lesssim \left(\|\hat{u}_0\|_{L^2(\Omega)}+\|u_0-\hat{u}_0\|_{H^3(\Omega)}\right)e^{-\kappa t}\\+\|u_0-\hat{u}_0\|_{L^2(\Omega)}e^{-\alpha t}.\end{multline}

\section{Numerics}
\subsection{Algorithm} In this section, we describe the steps to obtain the numerical solution of the plant-observer-error system given in \eqref{KdVBurgers}, \eqref{observer}, and \eqref{error}.  We follow a different approach compared to for instance \cite{Marx18}.  Our idea is based on first solving the models \eqref{error} and \eqref{hatw} with homogeneous boundary conditions and then obtaining the solutions of nonhomogeneous boundary value problems \eqref{KdVBurgers} and \eqref{observer} by using the invertibility of the backstepping transformation given in Lemma \ref{inverselem}.
\begin{description}
  \item[(Step 1)] At first we obtain numerical solutions of kernel models \eqref{kEq} and \eqref{p}. This is done via successive approximation. More precisely, we first change variables by setting $t\equiv y$, $s\equiv x-y$, and $G(s,t)\equiv k(x,y)$.  Then, $G$ satisfies the boundary value problem given by
\begin{eqnarray}
  \label{Geq}G_{ttt} - 3G_{stt}+ 3G_{sst} + G_{t} &=& -\lambda G, \\
  \label{Geqb}G(s,0) = G(0,t) &=& 0, \\
  \label{Geqc}G_s(0,t) &=& \frac{\lambda}{3}t
\end{eqnarray} on the triangular domain $\mathcal{T}_{0}\equiv \left\{ (s,t) \,|\, t \in [0,L], s \in [0,L-t]\right\}.$ Note that the solution of \eqref{Geq}-\eqref{Geqc} can be constructed by solving the integral equation
\begin{multline}\label{GepsInt}
  {G}(s,t) = \frac{{\lambda}}{3}st\\
  +\frac{1}{3}\int_0^t\int_0^s\int_0^\omega (-G_{ttt} + 3G_{stt}  - G_{t} -\lambda G)(\xi,\eta)d\xi d\omega d\eta.
\end{multline}  Therefore, we set \begin{multline}G^n(s,t)= \frac{{\lambda}}{3}st\\
+\frac{1}{3}\int_0^t\int_0^s\int_0^\omega (-G^{n-1}_{ttt} + 3G^{n-1}_{stt}  - G^{n-1}_{t} -\lambda G^{n-1})(\xi,\eta)d\xi d\omega d\eta\end{multline} for $n\ge 1$ with $G^0\equiv 0.$ We have proven in  \cite{BatalOzsari2018-1} that the sequence $G^n$ uniformly converges to a smooth function on $\mathcal{T}_{0}$. For the sake of numerical experiments, we define a parameter $n_{iter}\in \mathbb{Z}_+$ and use $$k_{num}(x,y)=G^{n_{iter}}(x-y,y)$$ for the kernel $k$. Since the solution of \eqref{p} is given by $p(x,y)=k(L-y,L-x)$, we will use $$p_{num}(x,y)=k_{num}(L-y,L-x)=G^{n_{iter}}(x-y,L-x)$$ for the kernel $p$. The observer gains $P_1$ and $P_2$ will then be taken as $$P_{1,num}(x)=\frac{\partial}{\partial_y}p_{num}(x,0)\text{ and }P_{2,num}(x,0)=-p_{num}(x,0).$$ Using these polynomial approximations, we also define approximations for $\Psi_i$, $i=1,2$ by setting $$\Psi_{i,num}(x)\equiv P_{i,num}(x)-\int_0^x P_{i,num}(y)k_{num}(x,y)dy.$$
\item[(Step 2)] Secondly, we numerically solve the error system \eqref{error}. In order to do this, we modify the finite difference scheme given in \cite{Pazato}.  To this end, we set the discrete space $$X_J:=\{\tilde{u}=(\tilde{u}_0,\tilde{u}_1,...,\tilde{u}_J)\in \mathbb{R}^{J+1}\,|\,\tilde{u}_0=\tilde{u}_{J-1}=\tilde{u}_J=0\},$$ and the difference operators $\displaystyle (D^+\tilde{u})_j:=\frac{\tilde{u}_{j+1}-\tilde{u}_j}{\delta x}$, $\displaystyle (D^-\tilde{u})_j:=\frac{\tilde{u}_{j}-\tilde{u}_{j-1}}{\delta x}$ for $j=1,...,J-1$, and $\displaystyle D=\frac{1}{2}(D^++D^-)$. Let $\delta x$ and $\delta t$ be the space and time steps for $j=0,...,J,$ and $n=0,1,...,N$, respectively. Then the numerical approximation of the linearised error system \eqref{error} takes the form
\begin{eqnarray}
  \label{wjn1}\frac{\tilde{u}_{j}^{n+1}-\tilde{u}_j^n}{\delta t}+(\mathcal{A}\tilde{u}^{n+1})_j= P_{1,num}(x_j)\frac{\tilde{u}_{1}^{n}}{\delta x}\\
  +P_{2,num}(x_j)\frac{(\tilde{u}_{2}^{n}-2\tilde{u}_{1}^{n})}{(\delta x)^2},\hspace{.1in} j=1,...,J-1\\
  \label{wjn2}\tilde{u}_0=\tilde{u}_{J-1}=\tilde{u}_J &=& 0, \\
  \label{wjn3} \tilde{u}_0  =\int_{x_{j-\frac{1}{2}}}^{x_{j^+\frac{1}{2}}}\tilde{u}_0(x)dx,\hspace{.1in} j=1,...,J-1,
\end{eqnarray} where $x_{j\mp\frac{1}{2}}=(j\mp\frac{1}{2})\delta x$, $x_j=j\delta x$. The $(J-1)\times (J-1)$ matrix $\mathcal{A}$ approximates $\tilde{u}_x+\tilde{u}_{xxx}$ and it is defined by $\mathcal{A}:=D^+D^+D^-+D$. Let us set $\tilde{\mathcal{C}}:=I+\delta t A$.  Then, from the main equation, we obtain
\begin{multline}\tilde{u}_{j}^{n+1}=\tilde{\mathcal{C}}^{-1}\left(\tilde{u}_j^n\right.\\\left.+P_{1,num}(x_j)\frac{(\delta t)\tilde{u}_{1}^{n}}{\delta x}+P_{2,num}(x_j)\frac{\delta t(\tilde{u}_{2}^{n}-2\tilde{u}_{1}^{n})}{(\delta x)^2}\right)\end{multline} for $j=1,...,J-1$.
\item[(Step 3)] The next step is to solve \eqref{hatw}.  The right hand side of the main equation in \eqref{hatw} includes the traces $\tilde{w}_x(0,t)$ and $\tilde{w}_{xx}(0,t)$.  Observe that these traces are equal to $\tilde{u}_x(0,t)$ and $\tilde{u}_{xx}(0,t)$ by the transformation \eqref{transtildew} and the boundary conditions $p(x,x)=0$ and $\tilde{w}(0,t)=0.$ Therefore, we can use the approximations $\displaystyle\frac{\tilde{u}_{1}^{n}}{\delta x}$ and $\displaystyle\frac{\tilde{u}_{2}^{n}-2\tilde{u}_{1}^{n}}{(\delta x)^2}$ from the previous step to approximate $\tilde{u}_x(0,t_n)$ and $\tilde{u}_{xx}(0,t_n)$ at the $n$\textsuperscript{th} time step.  Then the numerical approximation of the linearised observer target system \eqref{hatw} takes the form
\begin{eqnarray}
  \label{wjn1}\frac{\hat{w}_{j}^{n+1}-\hat{w}_j^n}{\delta t}+(\mathcal{A}\hat{w}^{n+1})_j+\lambda\hat{w}^{n+1}= (RHS),\\
  \label{wjn2}\hat{w}_0=\hat{w}_{J-1}=\hat{w}_J = 0, \\
  \label{wjn3} \hat{w}_0  =\int_{x_{j-\frac{1}{2}}}^{x_{j^+\frac{1}{2}}}\hat{w}_0(x)dx,
\end{eqnarray} for $j={1,...,J}$, where $\hat{w}_0$ is obtained from the transformation \eqref{transhatw} and \begin{multline}(RHS)=\frac{\partial}{\partial_y}k_{num}(x_j,0)-\Psi_{1,num}(x_j)\frac{\tilde{u}_{1}^{n}}{\delta x}\\-\Psi_{2,num}(x_j)\frac{(\tilde{u}_{2}^{n}-2\tilde{u}_{1}^{n})}{(\delta x)^2}.\end{multline} Let us set $\hat{\mathcal{C}}:=(1+\delta t\lambda)I+\delta t A$.  Then, from the main equation, we obtain
\begin{multline}\hat{w}_{j}^{n+1}=\hat{\mathcal{C}}^{-1}\left(\hat{w}_j^n+\frac{\partial}{\partial_y}k_{num}(x_j,0)-\Psi_{1,num}(x_j)\frac{\tilde{u}_{1}^{n}}{\delta x}\right.\\\left.-\Psi_{2,num}(x_j)\frac{(\tilde{u}_{2}^{n}-2\tilde{u}_{1}^{n})}{(\delta x)^2}\right)\end{multline} for $j=1,...,J-1$.

In order to obtain the solution of the observer system \eqref{observer}, we use the inverse of the transformation \eqref{transhatw}.  Given $\tilde{w}$, we can find the corresponding inverse image $\hat{u}$ via the succession method given in the proof of Lemma \ref{inverselem}  (see for example \cite[Lemma 2.4]{Liu03} and \cite[Lemma 2.2]{BatalOzsari2018-1}).  To this end, let $m_{iter}$ denote the number of iterations in the succession and set $v^0=\mathcal{K}\tilde{w}$, $v^{k}:=\mathcal{K}(\tilde{w}+v^{k-1})$ for $1\le k\le m_{iter}$, where $\mathcal{K}$ is the numerical approximation of the integral in the definition of $\Upsilon_k$.  Then, $v^{m_{iter}}$ is an approximation of $v=\Phi(\tilde{w})$, and one gets an approximation of the solution of the observer system by setting $\hat{u}(x_j,t_n):=\hat{w}(x_j,t_n)+v^{m_{iter}}(x_j,t_n)$.

\item[(Step 4)] Finally, we solve the original plant \eqref{KdVBurgers} by setting $${u}(x_j,t_n):=\hat{u}(x_j,t_n)+\tilde{u}(x_j,t_n).$$
\end{description}
\subsection{Simulations}
In this section, we give two simulations for the linear model on a domain of critical length: (i) uncontrolled solution and (ii) controlled solution.  The first simulation (Fig. 1) shows a time independent solution of the KdV equation on $\Omega=(0,2\pi)$ with initial datum $u_0=1-\cos x$ when no boundary feedback is present. This is the case when all boundary conditions are homogeneous: $u(0,t)=u(2\pi,t)=u_x(2\pi,t)=0$.  The second simulation (Fig. 2) shows the solution of the KdV equation with the same initial datum but subject to the backstepping feedback controllers given in \eqref{controllers} which use the state of the observer system.  The bump at $x=2\pi$ in Fig. 4 represents the action of the feedbacks at the right endpoint of the domain.
\begin{figure}
  \centering
   \includegraphics[scale=.60]{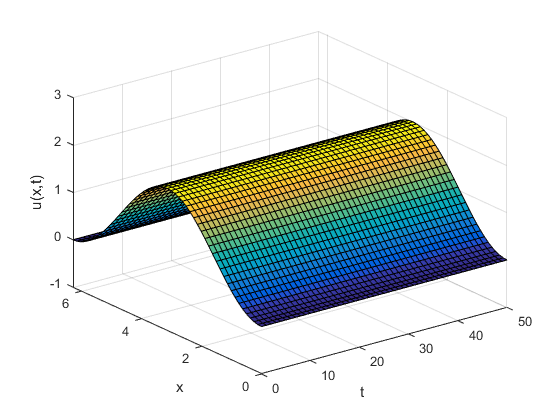}
  \caption{Uncontrolled solution with initial datum $u_0=1-\cos(x)$ on a domain of length $2\pi$.}\label{uncont-sol}
\end{figure}

\begin{figure}
  \centering
   \includegraphics[scale=.6]{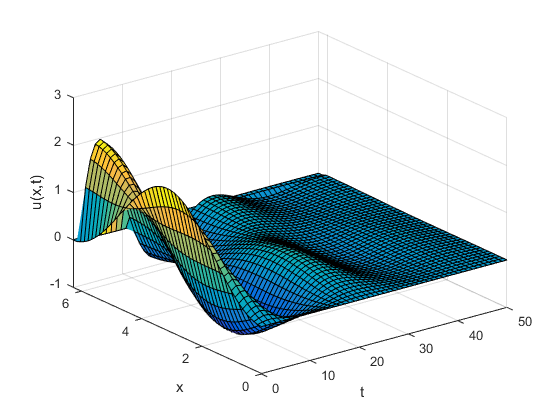}
  \caption{Controlled solution with initial datum $u_0=1-\cos(x),\hat{u}_0=0$, kernel parameters $\lambda=\tilde{\lambda}=0.01$ on a domain of length $L=2\pi$, $n_{iter}=m_{iter}=10$.}\label{uncont-sol}
\end{figure}

\section{Output feedback stabilization with a single controller and boundary measurement}
Using two feedback controllers at the right endpoint of the domain and measuring two traces at the left are not necessary to obtain the stabilization results in Section \ref{ObserverDes}. One can achieve this by using only one controller and making only one measurement as well.  More precisely, if we respectively take $V(t)=0$ and $P_1(x)=0$ in \eqref{KdVBurgers} and \eqref{observer}, then Theorem \ref{mainthm} still holds but with decay rate constants smaller than $\alpha$ and $\kappa$.  To see this, let us assume $u$ solves \eqref{KdVBurgers} with  $V(t)=0$, and $\hat{u}$ solves \eqref{observer} with  $V(t)=P_1(x)=0$. Then the error target system becomes
\begin{equation}\label{tildewnew} \left\{ \begin{array}{ll}
        \tilde{w}_t+\tilde{w}_x+\tilde{w}_{xxx}+\tilde{\lambda} \tilde{w}=-p_y(x,0)\tilde{w}_x(0,t), \text { in } \Omega\times (0,T),\\
        \tilde{w}(0,t)=0,\,\tilde{w}(L,t)=0,\,\\\tilde{w}_x(L,t)=\int_0^Lp_x(L,y)\tilde{w}(y,t)dy, \text { in } (0,T),\\
        \tilde{w}(x,0)=\tilde{w}_0(x), \text { in } \Omega. \end{array} \right.
\end{equation}
Applying the same multipliers to \eqref{tildewnew} as in Section \ref{ObserverDes}, we get
\begin{eqnarray}
\label{aaaa}\|\tilde{w}(t)\|_{L^2(\Omega)}\leq \|\tilde{w}(0)\|_{L^2(\Omega)} e^{-\beta t}, \\
\label{bbbb}\|\tilde{w}_t(t)\|_{L^2(\Omega)}\leq \|\tilde{w}_t(0)\|_{L^2(\Omega)} e^{-\beta t},
\end{eqnarray}
where $\beta=\big(\tilde{\lambda}-\frac{1}{2}\|p_x(L,\cdot )\|^2_{L^2(\Omega)}-\frac{1}{2}\|p_y(\cdot,0 )\|^2_{L^2(\Omega)}\big)$.
Moreover,
 $\|\tilde{w}_t\|_{L^2(\Omega)}\leq\|\tilde{w}_x+\tilde{w}_{xxx}+\tilde{\lambda}\tilde{w}\|_{L^2(\Omega)}$ $+
 \|p_y(\cdot, 0)\|_{L^2(\Omega)}|\tilde{w}_x(0,t)|;$
and $\tilde{w}_x(0,t)=\tilde{w}_x(L,t)-\int_0^L \tilde{w}_{xx}(x,t)dx,$
which, together with the boundary condition, implies
\begin{equation}
\label{anew2.5}
|\tilde{w}_x(0,t)|\lesssim \|\tilde{w}\|_{L^2(\Omega)}+\|\tilde{w}_{xx}\|_{L^2(\Omega)}.
\end{equation}
Therefore
\begin{equation}
\label{anew3}
\|\tilde{w}_t\|_{L^2(\Omega)}\lesssim \|\tilde{w}\|_{H^3(\Omega)}.
\end{equation}
Combining \eqref{bbbb} and \eqref{anew3}, we obtain
\begin{equation}
\label{anew4}
\|\tilde{w}_t(t)\|_{L^2(\Omega)}\leq \|\tilde{w}_0\|_{H^3(\Omega)} e^{-\beta t}.
\end{equation}
On the other hand by \eqref{tildewnew} and \eqref{anew2.5} we have
\begin{multline}\label{anew5}\|\tilde{w}_{xxx}\|^2_{L^2(\Omega)}\leq 4(\|\tilde{w}_x(t)\|^2_{L^2(\Omega)}+\tilde{\lambda}\|\tilde{w}(t)\|^2_{L^2(\Omega)}\\
+\|\tilde{w}_t(t)\|^2_{L^2(\Omega)}+\|p_y(\cdot, 0)\|^2_{L^2(\Omega)}|\tilde{w}_x(0,t)|^2)\\
\lesssim \|\tilde{w}_{x}(t)\|^2_{L^2(\Omega)}+\|\tilde{w}_{xx}(t)\|^2_{L^2(\Omega)}\\+\|\tilde{w}(t)\|^2_{L^2(\Omega)}+\|\tilde{w}_t(t)\|^2_{L^2(\Omega)}.
\end{multline}
Applying $\epsilon$-Young's inequlity to the  Gagliardo-Nirenberg inequality $\|\tilde{w}_{xx}\|_{L^2(\Omega)}\leq \|\tilde{w}_{xxx}\|^{\frac{2}{3}}_{L^2(\Omega)} \|\tilde{w}\|^{\frac{1}{3}}_{L^2(\Omega)}$, we obtain
\begin{equation}
\label{anew6}
\|\tilde{w}_{xx}\|^2_{L^2(\Omega)}\leq \delta\|\tilde{w}_{xxx}\|^{2}_{L^2(\Omega)}+ c_{\delta}\|\tilde{w}\|^{2}_{L^2(\Omega)}
\end{equation}
for any $\delta>0$. Combining \eqref{wxest}, \eqref{anew5} and \eqref{anew6}, we see that our new error target $\tilde{w}$ also satisfies \eqref{anewold1}. Hence \eqref{anewold2}, which together with \eqref{aaaa} and \eqref{anew4} implies
\begin{equation}\label{part2new}
\|\tilde{w}\|_{H^3(\Omega)} \lesssim \|\tilde{w}_0\|_{H^3(\Omega)}e^{-\beta t}.
\end{equation}
Not only $\|\tilde{w}\|_{H^3(\Omega)}$ but also $|\tilde{w}_{xx}(0,t)|^2$ is bounded by $\|\tilde{w}_0\|_{H^3(\Omega)}e^{-\beta t}$. To see this let us multiply \eqref{tildewnew} by $(L-x)\tilde{w}_{xx}$ and integrate over $\Omega$. After applying integration by parts and the boundary conditions we obtain
\begin{multline}
\tilde{w}_x^2(0,t)+\tilde{w}_{xx}^2(0,t)\\
=\frac{2}{L}\int_0^L \big( (L-x)\tilde{w}_t\tilde{w}_{xx}+\frac{1}{2}\tilde{w}_{x}^2+\frac{1}{2}\tilde{w}_{xx}^2+\tilde{\lambda}(L-x)\tilde{w}\tilde{w}_{xx}\big) dx\\
-\frac{2}{L}\tilde{w}_x(0,t)\int_0^L (L-x)p_y(x,0)\tilde{w}_{xx}dx.
\end{multline}
Using Cauchy-Schwarz and Young's inequalities, we achieve $
|\tilde{w}_{xx}(0,t)|^2 \lesssim \|\tilde{w}_t\|_{L^2(\Omega)}^2+\|\tilde{w}\|_{H^3(\Omega)}^2,$
which, together with \eqref{anew4} and \eqref{part2new}, implies
\begin{equation}\label{anew7}
|\tilde{w}_{xx}(0,t)|^2 \lesssim \|\tilde{w}_0\|_{H^3(\Omega)}e^{-\beta t}.
\end{equation}

In the case of one observer, i.e., $P_1(x)=V(t)=0$, $\hat{w}$ given by \eqref{transhatw} solves the following observer target system
\begin{equation}\label{hatwnew} \left\{ \begin{array}{ll}
        \hat{w}_t+\hat{w}_x+\hat{w}_{xxx}+\lambda \hat{w}\\
        =k_y(x,0)\hat{w}_x(0,t))-\Psi_2(x)\tilde{w}_{xx}(0,t),\\
        \hat{w}(0,t)=0, \; \hat{w}(L,t)=0, \;\\ \hat{w}_x(L,t)=-\int_0^L k_x(L,y)\hat{u}(y,t)dy. \end{array} \right.
\end{equation}
 Multiplying \eqref{hatwnew} by $\hat{w}$ and integrating over $\Omega$, we obtain
\begin{multline}
\label{anew8}
\frac{1}{2}\frac{d}{dt}\|\hat{w}(t)\|_{L^2(\Omega)}^2+\lambda\|\hat{w}(t)\|_{L^2(\Omega)}^2+\frac{1}{2}
|\hat{w}_x(0,t)|^2\\
=\hat{w}_x(0,t)\int_0^L k_y(x,0)\hat{w}(x,t)dx\\
-\tilde{w}_{xx}(0,t)\int_0^L\Psi_2(x) \hat{w}(x,t)dx\\+\frac{1}{2}\left[\int_0^L k_x(L,y)\hat{u}(y,t)dx\right]^2.
\end{multline}
Note that $\hat{w}=(I-\Upsilon_{k})\hat{u}$. Therefore by Lemma \ref{inverselem} we have $\|\hat{u}\|_{L^2(\Omega)}\leq \|(I-\Upsilon_{k})^{-1}\|_{B[L^2(\Omega)]}\|\hat{w}\|_{L^2(\Omega)}$. Using this fact and applying $\epsilon$-Young's and Cauchy-Schwarz inequalities to the right hand side of \eqref{anew8}, for any $\epsilon$ we get
\begin{equation}\label{estenergynew}
\frac{1}{2}\frac{d}{dt}\|\hat{w}(t)\|_{L^2(\Omega)}^2+\mu \|\hat{w}(t)\|_{L^2(\Omega)}^2
\leq \frac{1}{2\epsilon}\tilde{w}^2_{xx}(0,t),
\end{equation}
where \begin{multline*}
\mu\equiv \lambda-\frac{1}{2}\|k_y(\cdot,0)\|_{L^2(\Omega)}^2-\frac{1}{2}\epsilon\|\Psi_2\|_{L^2(\Omega)}^2\\
-\frac{1}{2}\|k_x(L,\cdot)\|_{L^2(\Omega)}^2\|(I-\Upsilon_{k})^{-1}\|_{B[L^2(\Omega)]}^2.
\end{multline*}
Inequality \eqref{estenergynew} and \eqref{part2new} imply
\begin{equation}\label{estenergynew2}
\frac{1}{2}\frac{d}{dt}\|\hat{w}(t)\|_{L^2(\Omega)}^2+\mu \|\hat{w}(t)\|_{L^2(\Omega)}^2
\lesssim \|\tilde{w}_0\|^2_{H^3(\Omega)}e^{-2\beta t}.
\end{equation}
By the proof of \cite[Lemma 2.5]{BatalOzsari2018-1}, we know that asymptotically $\|k_y(\cdot,0)\|_{L^2(\Omega)}\sim \lambda$. A similar argument also implies $\|k_x(L,\cdot)\|_{L^2(\Omega)}\sim \lambda$. Moreover, using the calculations in \cite{Liu03}, it is not hard to see that $\|(I-\Upsilon_{k})^{-1}\|_{B[L^2(\Omega)]} \sim 1+\lambda e^{C\lambda} $ where $C > 0$
depends only on $L$. Therefore choosing $\lambda$ and $\epsilon$ sufficiently small we can guarantee that $\mu>0$. In addition, in the case of $\lambda=\tilde{\lambda}$, we have $p(x,y)=k(L-y,L-x)$ and  $\|k_y(\cdot,0)\|_{L^2(\Omega)}=\|p_x(L,\cdot)\|_{L^2(\Omega)}$, $\|p_y(\cdot,0)\|_{L^2(\Omega)}=\|k_x(L,\cdot)\|_{L^2(\Omega)}$ which imply that choosing $\lambda=\tilde{\lambda}$ if necessary we can also guarantee $\beta>\mu.$ Taking $\beta>\mu$, multiplying \eqref{estenergynew2} by $e^{2\mu t}$ and integrating from $0$ to $t$ we obtain
\begin{equation}\label{expdecay1new}
\|\hat{w}(t)\|_{L^2(\Omega)}\lesssim \left(\|\hat{w}_0\|_{L^2(\Omega)}+\|\tilde{w}_0\|_{H^3(\Omega)}\right) e^{-\mu t},
\end{equation}
which together with \eqref{normeq} implies
\begin{equation}
\label{anew9}
\|\hat{u}\|_{L^2(\Omega)}\lesssim \left(\|\hat{u}_0\|_{L^2(\Omega)}+\|u_0-\hat{u}_0\|_{H^3(\Omega)}\right)e^{-\mu t}.
\end{equation}
By  \eqref{aaaa}, \eqref{part2new} and \eqref{normeq} we also have
\begin{eqnarray}\label{utilde2}
\|u-\hat{u}\|_{L^2(\Omega)} &\lesssim& \|u_0-\hat{u}_0\|_{L^2(\Omega)}e^{-\beta t},\\
\|u-\hat{u}\|_{H^3(\Omega)} &\lesssim& \|u_0-\hat{u}_0\|_{H^3(\Omega)}e^{-\beta t}.
\end{eqnarray}
Again, combining \eqref{anew9} and \eqref{utilde2} and using the triangle inequality, we prove the exponential decay of $u$.

\section{Conclusion}
In this paper, we studied an output feedback stabilization problem with right endpoint controller(s) to which the standard backstepping method does not apply because the associated kernel PDE models become overdetermined and do not possess smooth solutions. The difficulty was due to the type of given boundary conditions (one b.c. at the left, two b.c. at the right) and the location of the controller(s).  We dealt with this issue by using a kernel instead, that does not satisfy all of the boundary conditions implied by the standard algorithm of backstepping.  Although using such a kernel is associated with more complicated target systems and slower rate of decay, it had the major advantage that the exponential stabilization can be achieved even on critical length domains.  This method is interesting in the sense that it can be applied to many other PDEs where one encounters overdetermined kernel models.
\section*{Acknowledgement}
We would like to thank the anonymous reviewers for their careful reading of the manuscript and their several insightful comments and suggestions,
which significantly contributed to improving the quality of this article.
\bibliographystyle{plain}        
\bibliography{myreferences}           

\end{document}